\documentclass[12pt,a4paper,reqno]{amsart}
\usepackage{amsfonts,amsmath,amssymb}

\advance\textheight by \topskip
\newtheorem{theorem}{Theorem}
\newtheorem{conjecture}[theorem]{Conjecture}

\numberwithin{equation}{section}

\newcommand{\gauss}[2]{\genfrac{[}{]}{0pt}{}{#1}{#2}_q}

\title[New  $q$--tangent and $q$--secant numbers]
{Combinatorics of geometrically distributed random
variables:\\  New   $q$--tangent and $q$--secant numbers}

\author{Helmut Prodinger}

\thanks{Parts of this research were conducted while the author was 
a visitor of the Technical University of Graz where he was supported
by the start project Y96--MAT}

\address{ Helmut Prodinger,
Centre for Applicable Analysis and Number Theory,
 Department of Mathematics,
University of the Witwatersrand, P.~O. Wits, 
2050 Johannesburg, South Africa,  email:
{\tt helmut@gauss.cam.wits.ac.za},\newline
homepage: {\tt http://www.wits.ac.za/helmut/index.htm}
}

\date{October 18, 1999}

\begin{document}

\begin{abstract} Up--down permutations are counted by tangent resp.
secant numbers. Considering words instead, where the letters are
produced by independent geometric distributions, there are several
ways of introducing this concept; in the limit they all coincide with
the classical version. In this way, we get some new $q$--tangent
and $q$--secant functions. Some of them also have nice continued
fraction expansions; in one particular case,  we could not find a proof
for it. Divisibility results \`a la Andrews/Foata/Gessel are also
discussed. 
\end{abstract}

\maketitle

\section{Introduction}

Permutations $\pi=\pi_1\pi_2\cdots\pi_n$ are
called up--down permutations if
$\pi_1<\pi_2>\pi_3<\pi_4>\pi_5\cdots$.
For odd $n$, the number of them is given by
$n![z^n]\tan z$, and for even $n$ by
$n![z^n]\sec z$. 
One finds that in many textbooks, e.~g. \cite{GrKnPa94}.

Instead of speaking about exponential generating functions, we
prefer to think of the coefficients of
$\tan z$ and $\sec z$ as probabilities.

If we consider words $a_1a_2\dots a_n$ with letters in $\{1,2,\dots\}$ 
with probabilities (weights) $p, pq, pq^2,\dots$, where $p+q=1$
(independent geometric probilities), then there are several ways to
introduce this concept. We can use $<$ or $\le$ for ``up,''
$>$ or $\ge $ for ``down,'' which gives 4 possibilities. Also,
it {\it makes\/} a difference to consider ``up--down'' versus
``down--up.'' That gives in principle 8 versions for $q$--tangent
and $q$--secant numbers. However, reading the word from right to left, 
the instance ``$\le>\le>\dots$'' coincides with the ``$<\ge<\ge<\dots$,'' 
and similarly for ``$\ge<\ge<\dots$'' and ``$>\le>\le\dots$,'' which gives
us 6  $q$--tangent numbers (probabilities, to be more precise).
In the instance of even length (secant numbers), there are more symmetries, 
and we have only 4 $q$--secant numbers. 

By general principles, the limit $q\to1$ reduces all the instances to the
classical quantities.

We need a few definitions from $q$--analysis; consult the books
\cite{Andrews76} and \cite{AnAsRo99}:

\begin{gather}\begin{split}
[n]_q:=\frac{1-q^n}{1-q},\qquad [n]_q!:=[1]_q[2]_q\dots[n]_q,\\
(x;q)_n:=(1-x)(1-xq)(1-xq^2)\dots(1-xq^{n-1}).
\end{split}\end{gather}

\section{Recursions}

We introduce the functions

\begin{equation}
T_n^{\le\ge}(u)
\end{equation}
where the coefficient of $u^i$ in it is the probability that a word
of length $n$ satisfies the $\le\ge\le\ge\dots$ condition and
ends with the letter $i$. Also, we define

\begin{equation}
\tau_n^{\le\ge}= T_n^{\le\ge}(1),
\end{equation}
which drops the technical condition about the last letter. 

Furthermore, we introduce the generating functions

\begin{equation}
F^{\le\ge}(z,u)=\sum_{n\ge0}
T_n^{\le\ge}(u)z^n\qquad\text{and}\qquad
f^{\le\ge}(z)=F^{\le\ge}(z,1).
\end{equation}

Quantities like $F^{\le>}(z,u)$ etc. are defined in an obvious way.

For the instance of secant numbers, we define similar quantities, but
use the letters $S,\sigma,G,g$ instead of $T,\tau, F,f$.

Obviously we only get nonzero contributions for odd $n$ in the tangent
case and for even $n$ in the secant case. 

The reason to operate with a variable $u$ that controls the last letter
is the technique of ``adding a new slice,'' that was applied with success in 
\cite{FlPr87} and, more recently, in \cite{KnPr98}.

\begin{theorem}
The functions $T_{2n+1}^{\nabla\triangle}(u)$
satisfy the following recurrences:

\begin{itemize}

\item 
\begin{align} \begin{split}  
T^{\ge\le}_{2n+1}(u)&=
\frac{p^2u}{(1-qu)(1-q^2u)}
T^{\ge\le}_{2n-1}(1)
-\frac{p^2u}{(1-qu)(1-q^2u)}
T^{\ge\le}_{2n-1}(q^2u)\\
 T^{\ge\le}_{1}(u)&=
\frac{pu}{1-qu},
\end{split}\end{align}

\item

\begin{align}\begin{split}  
T^{\ge<}_{2n+1}(u)&=
\frac{p^2qu^2}{(1-qu)(1-q^2u)}
T^{\ge<}_{2n-1}(1)-
\frac{p^2qu^2}{(1-qu)(1-q^2u)}
T^{\ge<}_{2n-1}(q^2u)\\
T^{\ge<}_{1}(u)&=
\frac{pu}{1-qu}.
\end{split}\end{align}

\item

\begin{align}\begin{split}  
T^{><}_{2n+1}(u)&=
\frac{p^2qu^2}{(1-qu)(1-q^2u)}
T^{><}_{2n-1}(1)-
\frac{p^2u}{q(1-qu)(1-q^2u)}
T^{><}_{2n-1}(q^2u)\\
T^{><}_{1}(u)&=
\frac{pu}{1-qu},
\end{split}\end{align}

\item
\begin{align}\begin{split}  
T^{\le\ge}_{2n+1}(u)&=
\frac{pu}{q(1-qu)}
T^{\le\ge}_{2n-1}(q)-
\frac{p^2u}{q(1-qu)(1-q^2u)}
T^{\le\ge}_{2n-1}(q^2u)\\
 T^{\le\ge}_{1}(u)&=
\frac{pu}{1-qu},
\end{split}\end{align}

\item
\begin{align}\begin{split}  
T^{\le>}_{2n+1}(u)&=
\frac{pu}{q(1-qu)}
T^{\le>}_{2n-1}(q)-
\frac{p^2}{q^2(1-qu)(1-q^2u)}
T^{\le>}_{2n-1}(q^2u)\\
T^{\le>}_{1}(u)&=
\frac{pu}{1-qu},
\end{split}\end{align}

\item
\begin{align}\begin{split}  
T^{<>}_{2n+1}(u)&=
\frac{pu}{1-qu}
T^{<>}_{2n-1}(q)-
\frac{p^2u}{(1-qu)(1-q^2u)}
T^{<>}_{2n-1}(q^2u)\\
T^{<>}_{1}(u)&=
\frac{pu}{1-qu}.
\end{split}\end{align}

\end{itemize}

\end{theorem}

\begin{proof}
Since the technique is the same for all the instances, it is enough to discuss
e.~g. the ``$\ge\le$'' case. Adding a new slice means adding a pair $(k,j)$
with $1\le k\le i$, $j\ge k$, replacing $u^i$ by 1 and providing the factor
$u^j$. But
\begin{equation*}
\sum_{k=1}^ipq^{k-1}\sum_{j\ge k}pq^{j-1}u^j=
\frac{p^2u}{(1-qu)(1-q^2u)}-\frac{p^2u}{(1-qu)(1-q^2u)}\big(q^2u\big)^i,
\end{equation*}

which explains the recursion. The starting value is just
\begin{equation*}
\sum_{j\ge 1}pq^{j-1}u^j=\frac{pu}{1-qu}.
\end{equation*}

\end{proof}

\begin{theorem}\label{tan_theorem}The numbers  
 $\tau_{2n+1}^{\nabla\triangle}$
have the generating functions $f^{\nabla\triangle}(z)=\tan_q(z)
=\sin_q(z)/\cos_q(z)$:
\begin{itemize}
\item
\begin{equation}   
f^{\ge\le}(z)={
\sum\limits_{n\ge0}\dfrac{(-1)^nz^{2n+1}}{[2n+1]_q!}
q^{n(n+1)}
}\bigg/{
\sum\limits_{n\ge0}\dfrac{(-1)^nz^{2n}}{[2n]_q!}
q^{n(n-1)} 
}
\end{equation}
\item
\begin{equation}   
f^{\ge<}(z)={
\sum\limits_{n\ge0}\dfrac{(-1)^nz^{2n+1}}{[2n+1]_q!}
}\bigg/{
\sum\limits_{n\ge0}\dfrac{(-1)^nz^{2n}}{[2n]_q!}
}
\end{equation}
\item
\begin{equation}   
f^{><}(z)={
\sum\limits_{n\ge0}\dfrac{(-1)^nz^{2n+1}}{[2n+1]_q!}
q^{n^2}
}\bigg/{
\sum\limits_{n\ge0}\dfrac{(-1)^nz^{2n}}{[2n]_q!}
q^{n^2} 
}
\end{equation}

\item
\begin{equation}   
f^{\le\ge}(z)={
\sum\limits_{n\ge0}\dfrac{(-1)^nz^{2n+1}}{[2n+1]_q!}
q^{n^2}
}\bigg/{
\sum\limits_{n\ge0}\dfrac{(-1)^nz^{2n}}{[2n]_q!}
q^{n(n-1)} 
}
\end{equation}

\item
\begin{equation}   
f^{\le>}(z)={
\sum\limits_{n\ge0}\dfrac{(-1)^nz^{2n+1}}{[2n+1]_q!}
}\bigg/{
\sum\limits_{n\ge0}\dfrac{(-1)^nz^{2n}}{[2n]_q!}
}
\end{equation}
\item
\begin{equation}   
f^{<>}(z)={
\sum\limits_{n\ge0}\dfrac{(-1)^nz^{2n+1}}{[2n+1]_q!}
q^{n(n+1)}
}\bigg/{
\sum\limits_{n\ge0}\dfrac{(-1)^nz^{2n}}{[2n]_q!}
q^{n^2} 
}
\end{equation}

\end{itemize}
\end{theorem}

\begin{proof}
The proofs of the first 3 relations are very  similar, and we only
sketch the first instance. Summing up we find

\begin{equation*}
F^{\ge\le}(z,u)=
\frac{puz}{1-qu}+
\frac{p^2uz^2}{(1-qu)(1-q^2u)}
F^{\ge\le}(z,1)
-\frac{p^2uz^2}{(1-qu)(1-q^2u)}
F^{\ge\le}(z,q^2u)
\end{equation*}

Iterating that we find for $f(z)=f^{\ge\le}(z)$:

\begin{align*}
f(z)&=\frac{pz}{1-q}+\frac{p^2z^2}{(1-q)(1-q^2)}f(z)
-\frac{p^2q^2z^3}{(1-q)(1-q^2)(1-q^3)}\\&-
\frac{p^4q^2z^4}{(1-q)(1-q^2)(1-q^3)(1-q^4)}f(z)+\dots
\end{align*}

from which the announced formula follows by solving for $f(z)$.

The 3 others are trickier, because of a term $T_{2n-1}^{\nabla\Delta}(q)$.
Again, let us discuss one case. Observe that
\begin{equation*}
T_{2n-1}^{\le\ge}(q)=q^{-1}S_{2n}^{\le\ge}(1), 
\end{equation*}
because one more ``up'' step should replace $u^i$ by $\sum_{k\ge i}pq^{k-1}=
q^{i-1}$. Now the generating function $g^{\le\ge}(z)$ of the quantities $S_{2n}^{\le\ge}(1)$
(upcoming) is obtained independently, whence we get
\begin{equation*}
F^{\le\ge}(z,u)=\frac{puz}{q^2(1-qz)}g^{\le\ge}(z)
-\frac{p^2uz^2}{q(1-qu)(1-q^2u)}F^{\le\ge}(z,q^2u).
\end{equation*}
Now iteration as usual derives the desired result.
\end{proof}

\begin{theorem}
The functions $S_{2n}^{\nabla\triangle}(u)$
satisfy the following recurrences:

\begin{itemize}

\item 
\begin{align}   \begin{split}  
S^{\le\ge}_{2n+2}(u)&=
\frac{p^2u}{(1-qu)(1-q^2u)}
S^{\le\ge}_{2n}(1)-
\frac{p^2u}{(1-qu)(1-q^2u)}
S^{\le\ge}_{2n}(q^2u)\\
S^{\le\ge}_{2}(u)&=
\frac{p^2u}{(1-qu)(1-q^2u)},
\end{split}\end{align}

\item

\begin{align}\begin{split}  
S^{\le>}_{2n+2}(u)&=
\frac{p^2u}{(1-qu)(1-q^2u)}
S^{\le>}_{2n}(1)-
\frac{p^2}{q^2(1-qu)(1-q^2u)}
S^{\le>}_{2n}(q^2u)\\
S^{\le>}_{2}(u)&=
\frac{p^2u}{(1-qu)(1-q^2u)},
\end{split}\end{align}
\item

\begin{align}\begin{split}  
S^{<\ge}_{2n+2}(u)&=
\frac{p^2qu^2}{(1-qu)(1-q^2u)}
S^{<\ge}_{2n}(1)-
\frac{p^2qu^2}{(1-qu)(1-q^2u)}
S^{<\ge}_{2n}(q^2u)\\
S^{<\ge}_{2}(u)&=
\frac{p^2qu^2}{(1-qu)(1-q^2u)},
\end{split}\end{align}

\item

\begin{align}\begin{split}  
S^{<>}_{2n+2}(u)&=
\frac{p^2qu^2}{(1-qu)(1-q^2u)}
S^{<>}_{2n}(1)-
\frac{p^2u}{q(1-qu)(1-q^2u)}
S^{<>}_{2n}(q^2u)\\
S^{<>}_{2}(u)&=
\frac{p^2qu^2}{(1-qu)(1-q^2u)}.
\end{split}\end{align}

\end{itemize}

\end{theorem}
\begin{proof}
The proof works as in the easy cases of the tangent recursions and is
omitted. For the starting value, we must consider the first pair of
numbers. 
\end{proof}

\begin{theorem}The numbers  
 $\sigma_{2n}^{\nabla\triangle}$
have the generating functions $g^{\nabla\triangle}(z)=1/\cos_q(z)$:
\begin{itemize}
\item
\begin{equation}   
g^{\le\ge}(z)={
1}\bigg/{
\sum\limits_{n\ge0}\dfrac{(-1)^nz^{2n}}{[2n]_q!}
q^{n(n-1)} 
}
\end{equation}

\item
\begin{equation}   
g^{\le>}(z)={
1}\bigg/{
\sum\limits_{n\ge0}\dfrac{(-1)^nz^{2n}}{[2n]_q!}
			}
\end{equation}

\item
\begin{equation}   
g^{<\ge}(z)={
1}\bigg/{
\sum\limits_{n\ge0}\dfrac{(-1)^nz^{2n}}{[2n]_q!}
q^{n(2n-1)} 
}
\end{equation}

\item
\begin{equation}   
g^{<>}(z)={
1}\bigg/{
\sum\limits_{n\ge0}\dfrac{(-1)^nz^{2n}}{[2n]_q!}
q^{n^2} 
}
\end{equation}
												
\end{itemize}
\end{theorem}

\begin{proof}
The proofs are quite similar as before; however, iteration must be done for
the function $G^{\nabla\Delta}(z,u)-1$, and 1 must be added at the end. 
\end{proof}

\section{Jackson's $q$--sine and $q$--cosine functions}

Jackson in \cite{Jackson04} has introduced the functions

\begin{align}\begin{split}
\sin_q(z)&=\sum_{n\ge0}\frac{(-1)^nz^{2n+1}}{[2n+1]_q!},\\
\cos_q(z)&=\sum_{n\ge0}\frac{(-1)^nz^{2n}}{[2n]_q!}.
\end{split}
\end{align}
and proved the relation
\begin{equation}\label{sincos}
\sin_q(z)\sin_{1/q}(z)+\cos_q(z)\cos_{1/q}(z)=1
\end{equation}

Since we have here several $q$--sine and $q$--cosine functions, we call them
a $q$--sine--cosine pair, if relation (\ref{sincos}) holds.

\begin{theorem}\label{akin}
For the functions
$$
\sin_q(z):=\sum_{n\ge0}\frac{(-1)^nz^{2n+1}}{[2n+1]_q!}q^{An^2+Bn}
$$
and

$$
\cos_q(z):=\sum_{n\ge0}\frac{(-1)^nz^{2n}}{[2n]_q!}q^{Cn^2+Dn}
$$
exactly the  12 pairs in Table~\ref{tab1} are $q$--sine--cosine pairs:

\begin{table}[h]

\tabcolsep4mm
\arraycolsep4mm

\begin{center}
 
\begin{tabular}{  | c | c | c | c | }
  \hline
   \rule[-1.2ex]{0mm}{4ex}
   $A$   &   $B$    & $C$   &$D$   \\
  \hline
   \rule[-1.2ex]{0mm}{4ex}   
$0$&$0$&$0$&$0$\\
\hline
   \rule[-1.2ex]{0mm}{4ex}
$2$&$1$&$0$&$0$\\
\hline
   \rule[-1.2ex]{0mm}{4ex}
$0$&$0$&$2$&$-1$\\
\hline
   \rule[-1.2ex]{0mm}{4ex}
$2$&$1$&$2$&$-1$\\
\hline
   \rule[-1.2ex]{0mm}{4ex}
$0$&$1$&$0$&$1$\\
\hline
   \rule[-1.2ex]{0mm}{4ex}
$2$&$0$&$0$&$1$\\
\hline
   \rule[-1.2ex]{0mm}{4ex}
$1$&$0$&$1$&$0$\\
\hline
   \rule[-1.2ex]{0mm}{4ex}
$1$&$1$&$1$&$0$\\
\hline
   \rule[-1.2ex]{0mm}{4ex}
$0$&$1$&$2$&$-2$\\
\hline
   \rule[-1.2ex]{0mm}{4ex}
$2$&$0$&$2$&$-2$\\
\hline
   \rule[-1.2ex]{0mm}{4ex}
$1$&$0$&$1$&$-1$\\
\hline
   \rule[-1.2ex]{0mm}{4ex}
$1$&$1$&$1$&$-1$\\
     \hline
\end{tabular}   

\end{center}

\caption{\label{tab1}}

\end{table}

\end{theorem}

\begin{proof}
The desired relation gives us more and more restrictions when we look
at the coefficients of $z^{2n}$. By a tedious search that will not be
reported here we find these 12 possibilities, and   all others can be
excluded. The  proof that this indeed
works is very similar for all of them, so we give just one, namely the
instance $(1,0,1,0)$.

Note the following expansions:

\begin{align*}
\sin_{1/q} z&=\sum_{n\ge0}\frac{(-1)^nz^{2n+1}}{[2n+1]_q!}q^{\binom{2n+1}{2}-n^2},\\
\cos_{1/q} z&=\sum_{n\ge0}\frac{(-1)^nz^{2n}}{[2n]_q!}q^{\binom{2n}{2}-n^2}.
\end{align*}

So we must prove that for $n\ge1$

\begin{equation*}
\sum_{k=0}^n\gauss{2n}{2k}q^{\binom{2k}{2}-k^2+(n-k)^2}=
\sum_{k=0}^{n-1}\gauss{2n}{2k+1}q^{\binom{2k+1}{2}-k^2+(n-k-1)^2}
\end{equation*}

or, reversing the order of summation in the second sum,

\begin{equation*}
\sum_{k=0}^n\gauss{2n}{2k}q^{2k^2-2k-2nk}=
\sum_{k=0}^{n-1}\gauss{2n}{2k+1}q^{2k^2+k-n-2nk}.
\end{equation*}
We rewrite this again as

\begin{equation*}
\sum_{k\text{ even}}\gauss{2n}{k}q^{-nk+\binom k2}=
\sum_{k\text{ odd}}\gauss{2n}{k}q^{-nk+\binom k2}.
\end{equation*}

Therefore we have to prove that

\begin{equation*}
\sum_{k=0}^{2n}\gauss{2n}{k}(-1)^kq^{-nk+\binom k2}=0.
\end{equation*}

We use the formula (10.0.9) in \cite{AnAsRo99}, see also
\cite{Andrews76},
\begin{equation*}
\sum_{k=0}^n\gauss nk
z^kq^{\binom k2}=\prod_{j=0}^{n-1}(1+q^jz).
\end{equation*}
The desired result now follows by replacing $n$ by
$2n$ and  plugging in $z=-q^{-n}$.
\end{proof}

\begin{theorem} The  6 $\tan_q(z)$ functions in 
Theorem~\ref{tan_theorem} all involve $q$--sine--cosine pairs. 
\end{theorem}

{\bf Remark.} Replacing $q$ by $1/q$ in the  $q$--sine--cosine pairs
and rewriting everything again in the $q$--notation means replacing the
vector $(A,B,C,D)$ of exponents by$(2-A,1-B,2-C,-1-D)$.

\begin{table}[h]

\tabcolsep4mm
\arraycolsep4mm

\begin{center}
 
\begin{tabular}{  | c | c | c | c | c | c | c | c | }
  \hline
   \rule[-1.2ex]{0mm}{4ex}
   $A$   &   $B$    & $C$   &$D$ &$A'$   &   $B'$    & $C'$   &$D'$  \\
  \hline
   \rule[-1.2ex]{0mm}{4ex}   
0&0&0&0&2&1&2&-1\\
\hline
   \rule[-1.2ex]{0mm}{4ex}
2&1&0&0&0&0&2&-1\\
\hline
   \rule[-1.2ex]{0mm}{4ex}
0&1&0&1&2&0&2&-2\\
\hline
   \rule[-1.2ex]{0mm}{4ex}
2&0&0&1&0&1&2&-2\\
\hline
   \rule[-1.2ex]{0mm}{4ex}
1&0&1&0&1&1&1&-1\\
\hline
   \rule[-1.2ex]{0mm}{4ex}
1&1&1&0&1&0&1&-1\\
\hline
   \end{tabular}   

\end{center}

\caption{\label{tab2}}

\end{table}

This reduces the 12 pairs to 6 pairs.

\section{continued fractions}

Some of the 12 tangent functions have nice continued fraction expansions. 

\begin{theorem}
For $(A,B,C,D)=(0,0,0,0)$ and $(A,B,C,D)=(2,1,2,-1)$ we have

\begin{equation}
\tan_q(z)=\cfrac{z}{[1]_qq^{-0}-\cfrac{z^2}{[3]_qq^{-1}-
\cfrac{z^2}{[5]_qq^{-2}-\cfrac{z^2}{[7]_qq^{-3}-\cfrac{z^2}{\ddots}}}}}
\end{equation}

The two tangent functions coincide, which is classical, since Jackson
\cite{Jackson04} has shown that for his functions
\begin{equation*}
\sin_qz\cos_{1/q}z-\sin_{1/q}z\cos_{q}z=0
\end{equation*}
holds. 
\end{theorem}
\begin{proof}
For the proof by induction we must do the following:
Set $a_n=[2n-1]_qq^{1-n}$ and
\begin{align*}
p_n(z)&=a_np_{n-1}(z)-z^2p_{n-2}(z),\qquad p_0(z)=0,\  p_1(z)=z,\\
q_n(z)&=a_nq_{n-1}(z)-z^2q_{n-2}(z),\qquad q_0(z)=1, \  q_1(z)=a_1.
\end{align*}
We must show that 
\begin{equation*}
[z^k]\Big(p_n(z)\cos_q z-q_n(z)\sin_qz\Big)=0 \quad\text{for } k\le 2n.
\end{equation*}

Now look at
\begin{equation*}
[z^k]\Big((a_np_{n-1}(z)-z^2p_{n-2}(z))\cos_q z
-(a_nq_{n-1}(z)-z^2q_{n-2}(z))\sin_qz\Big).
\end{equation*}
By the induction hypothesis we only have to show that
\begin{equation*}
[z^{2n-1}]\Big(p_{n}(z)\cos_q z
-q_{n}(z)\sin_qz\Big)=0.
\end{equation*}

However, we can easily show by induction that
$$
p_n(z)=\sum_kz^{2k+1}(-1)^k\frac{[2n-2k-1]_q!q^{k(2k+1)-\binom n2}}{[n-2k-1]_q![2k+1]_q!\prod_{i=1}^{n-1-2k}(1+q^i)}
$$
and

$$
q_n(z)=\sum_kz^{2k}(-1)^k\frac{[2n-2k]_q!q^{k(2k-1)-\binom n2}}{[n-2k]_q![2k]_q!\prod_{i=1}^{n-2k}(1+q^i)}
$$
holds (the hard part is to find these formul\ae).
We have to prove that
\begin{equation*}
\sum_{k\ge0}[z^{2k+1}]p_n(z)[z^{2n-2k-2}]\cos_qz
=\sum_{k\ge0}[z^{2k}]q_n(z)[z^{2n-2k-1}]\sin_qz
\end{equation*}
or

\begin{multline*}
\sum_{k=0}^{\lfloor\frac{n-1}2\rfloor}
\frac{[2n-2k-1]_q!q^{k(2k+1)}}{[n-2k-1]_q![2k+1]_q!\prod_{i=1}^{n-1-2k}(1+q^i)[2n-2k-2]_q!}
\\=\sum_{k=0}^{\lfloor\frac{n}2\rfloor}
\frac{[2n-2k]_q!q^{k(2k-1)}}{[n-2k]_q![2k]_q!\prod_{i=1}^{n-2k}(1+q^i)
[2n-2k-1]_q!}.
\end{multline*}
Thus we must prove

\begin{multline*}
\sum_{k=0}^{\lfloor\frac{n-1}2\rfloor}
\frac{(1-q^{2n-2k-1})q^{k(2k+1)}}{[n-2k-1]_q![2k+1]_q!\prod_{i=1}^{n-1-2k}(1+q^i)}
=\sum_{k=0}^{\lfloor\frac{n}2\rfloor}
\frac{(1-q^{2n-2k})q^{k(2k-1)}}{[n-2k]_q![2k]_q!\prod_{i=1}^{n-2k}(1+q^i)}
\end{multline*}
or
\begin{equation*}
\sum_{k=0}^{n}
\frac{(1-q^{2n-k})q^{\binom k2}(-1)^k}{[n-k]_q![k]_q!\prod_{i=1}^{n-k}(1+q^i)}=0,
\end{equation*}
or
\begin{equation*}
\sum_{k=0}^{n}\frac{1}{(q;q)_k(q^2;q^2)_{n-k}}
(1-q^{2n-k})q^{\binom k2}(-1)^k=0.
\end{equation*}

Now
\begin{align*}
\sum_{k=0}^{n}\frac{1}{(q;q)_k(q^2;q^2)_{n-k}}
q^{\binom k2}(-1)^k&=[z^n]\sum_{k\ge0}\frac{q^{\binom k2}(-1)^kz^k}{(q;q)_k}
\sum_{k\ge0}\frac{z^k}{(q^2;q^2)_{k}}\\
&=[z^n]\prod_{k\ge 0}(1-zq^k)\bigg/\prod_{k\ge 0}(1-zq^{2k})\\
&=[z^n]\prod_{k\ge 0}(1-zq^{2k+1})\\
&=q^n[z^n]\prod_{k\ge 0}(1-zq^{2k})\\
&=\frac{q^nq^{2\binom n2}(-1)^n}{(q^2;q^2)_n}=\frac{(-1)^nq^{n^2}}{(q^2;q^2)_n}.
\end{align*}

Similarly,
\begin{align*}
\sum_{k=0}^{n}\frac{q^{2n-k}}{(q;q)_k(q^2;q^2)_{n-k}}
q^{\binom k2}(-1)^k&=q^{2n}[z^n]\sum_{k\ge0}\frac{q^{\binom k2}(-1)^k(z/q)^k}{(q;q)_k}
\sum_{k\ge0}\frac{z^k}{(q^2;q^2)_{k}}\\
&=q^{2n}[z^n]\prod_{k\ge 0}(1-zq^{k-1})\bigg/\prod_{k\ge 0}(1-zq^{2k})\\
&=q^{2n}[z^n]\prod_{k\ge 0}(1-zq^{2k-1})\\
&=q^{2n}q^{-n}[z^n]\prod_{k\ge 0}(1-zq^{2k})\\
&=
\frac{q^nq^{2\binom n2}(-1)^n}{(q^2;q^2)_n}=\frac{(-1)^nq^{n^2}}{(q^2;q^2)_n}.
\end{align*}
This finishes the proof.

The  continued fraction for $(2,1,2,-1)$ follows by replacing $q$ by $1/q$.
\end{proof}

\begin{theorem}
For $(A,B,C,D)=(0,1,0,1)$  we have

\begin{equation}
\tan_q(z)=\cfrac{z}{[1]_qq^{-0}-\cfrac{z^2}{[3]_qq^{-2}-
\cfrac{z^2}{[5]_qq^{-2}-\cfrac{z^2}{[7]_qq^{-4}-\cfrac{z^2}{\ddots}}}}}
\end{equation}

The negative powers of $q$ go like $0,2,2,4,4,6,6,8,8,\dots$.
\end{theorem}

\begin{proof}
The proof follows the same lines; this time the polynomials (continuants) are
$$
p_n(z)=\sum_kz^{2k+1}(-1)^k\frac{[2n-2k-1]_q!q^{2k(k+1)-\binom n2-\lfloor\frac n2\rfloor}}{[n-2k-1]_q![2k+1]_q!\prod_{i=1}^{n-1-2k}(1+q^i)}
$$
and

$$
q_n(z)=\sum_kz^{2k}(-1)^k\frac{[2n-2k]_q!q^{2k^2-\binom n2-\lfloor\frac n2\rfloor}}{[n-2k]_q![2k]_q!\prod_{i=1}^{n-2k}(1+q^i)}.
$$
Hence we have to prove that
\begin{multline*}
\sum_{k=0}^{\lfloor\frac{n-1}2\rfloor}
\frac{(1-q^{2n-2k-1})q^{2k(k+1)-k}}{[n-2k-1]_q![2k+1]_q!\prod_{i=1}^{n-1-2k}(1+q^i)}
=\sum_{k=0}^{\lfloor\frac{n}2\rfloor}
\frac{(1-q^{2n-2k})q^{2k^2-k}}{[n-2k]_q![2k]_q!\prod_{i=1}^{n-2k}(1+q^i)};
\end{multline*}
from here on we can use the previous proof.

An alternative proof is by noting that 
\begin{equation*}   
\tan_q^{(0,1,0,1)}(z)=\frac{1}{\sqrt q}\tan_q^{(0,0,0,0)}(z\sqrt q)
\end{equation*}
and using the previous result. 
\end{proof}

\begin{theorem}
For $(A,B,C,D)=(2,0,2,-2)$  we have

\begin{equation}
\tan_q(z)=\cfrac{z}{[1]_qq^{-0}-\cfrac{z^2}{[3]_qq^{-0}-
\cfrac{z^2}{[5]_qq^{-2}-\cfrac{z^2}{[7]_qq^{-2}-\cfrac{z^2}{\ddots}}}}}
\end{equation}

The negative powers of $q$ go like $0,0,2,2,4,4,6,6,8,8,\dots$.
\end{theorem}
\begin{proof}This follows from the previous theorem by replacing
$q$ by $1/q$.
\end{proof}

\begin{conjecture}
For $(A,B,C,D)=(1,0,1,0)$  we have

\begin{equation}
\tan_q(z)=\cfrac{z}{[1]_qq^{0}-\cfrac{z^2}{[3]_qq^{-2}-
\cfrac{z^2}{[5]_qq^{1}-\cfrac{z^2}{[7]_qq^{-9}-\cfrac{z^2}{\ddots}}}}}
\end{equation}

The positive powers of $q$ go like $0,1,6,15,\dots$ ($k(2k-1)$).

The negative powers of $q$ go like $2,9,20,35\dots$ ($(k+1)(2k-1)$).
\end{conjecture}

{\it Comment.\/}
It might be useful to rewrite the continued fraction as

\begin{equation}
\cfrac{z}{1-\cfrac{z^2b_1}{1-
\cfrac{z^2b_2}{1-\cfrac{z^2b_3}{1-\cfrac{z^2b_4
}{\ddots}}}}}
\end{equation}
with
\begin{align*}   
b_k&=\frac{1}{[k]_q[k+1]_q}q^{-k+(-1)^k(2k-1)}\\
&=\frac{1}{[k]_q[k+1]_q}\bigg[\frac12q^{-3k+1}\big(1+q^{4k-2}\big)-	
\frac{(-1)^k}2q^{-3k+1}\big(1-q^{4k-2}\big)\bigg].
\end{align*}

The recursions for the continuants are now
\begin{align*}
p_n(z)&=p_{n-1}(z)-b_{n-1}z^2p_{n-2}(z),\qquad p_0(z)=0,\  p_1(z)=z,\\
q_n(z)&=q_{n-1}(z)-b_{n-1}z^2q_{n-2}(z),\qquad q_0(z)=1,\  q_1(z)=1.
\end{align*}

Unfortunately, even with this form, I am currently unable to guess the
coefficients of these polynomials, whence I must leave this expansion
as an open problem.

\begin{conjecture}
For $(A,B,C,D)=(1,1,1,-1)$  we have

\begin{equation}
\tan_q(z)=\cfrac{z}{[1]_qq^{-0}-\cfrac{z^2}{[3]_qq^{0}-
\cfrac{z^2}{[5]_qq^{-5}-\cfrac{z^2}{[7]_qq^{3}-\cfrac{z^2}{\ddots}}}}}
\end{equation}

The positive powers of $q$ go like $0,3,10,21\dots$ ($(k-1)(2k-1)$).

The negative powers of $q$ go like $0,5,14,27\dots$ ($(k-1)(2k+1)$).
\end{conjecture}

{\it Comment.\/} This would be a corollary of the previous expansion.

{\bf Remark.} Normally, as e.~g. in \cite{Flajolet80} and
\cite{HaRaZe99}, the continued fraction expansions of the
{\it ordinary generating function\/} of the tangent and
secant numbers are considered, whereas we stick here to
the {\it exponent (or probability) generating functions}.

\section{Divisibility}

\begin{theorem}
The coefficient

\begin{equation*}
[2n+1]_q[z^{2n+1}]\tan_q(z)
\end{equation*}
is divisible by
\begin{equation*}
(1+q)(1+q^2)\dots(1+q^n)
\end{equation*}
for the vectors of exponents $(0,0,0,0)$, $(2,1,2,-1)$, $(0,1,0,1)$,
$(2,0,2,-2)$, $(1,0,1,0)$,
$(1,1,1,-1)$.

\end{theorem}

\begin{proof} 
The proof of \cite{AnGe78} covers the first 4 instances, since we note that
\begin{equation*}   
\tan_q^{(0,1,0,1)}(z)=\frac{1}{\sqrt q}\tan_q^{(0,0,0,0)}(z\sqrt q).
\end{equation*}
							
The only open case is thus $(1,0,1,0)$, as the remaining one would follow
from duality. 
Thus, let us now consider 

\begin{align*}
\sin_{q} z&=\sum_{n\ge0}\frac{(-1)^nz^{2n+1}}{[2n+1]_q!}q^{n^2},\\
\cos_{q} z&=\sum_{n\ge0}\frac{(-1)^nz^{2n}}{[2n]_q!}q^{n^2},\\
\end{align*}
and $\tan_qz=\sin_qz/\cos_qz$.

We need the following computation that is akin to the one in
Theorem~\ref{akin}.
\begin{align*}
[z^{2n+1}]\sin_{1/q}z\cos_qz&=\sum_{k=0}^n
\frac{(-1)^{k}q^{k(k+1)}}{[2k+1]_q!}\frac{(-1)^{n-k}q^{(n-k)^2}}{[2n-2k]_q!}\\
&=\frac{q^{n^2}(-1)^n}{[2n+1]_q!}\sum_{k=0}^n
\gauss{2n+1}{2k+1}q^{2k^2+k-2nk}\\
&=\frac{q^{n(n+1)}(-1)^n}{[2n+1]_q!}\sum_{k\text{ odd}}
\gauss{2n+1}{k}q^{\binom k2-nk}\\
&=\frac{q^{n(n+1)}{(-1)^n}}{[2n+1]_q!}\frac12\sum_{k=0}^n
\gauss{2n+1}{k}q^{\binom k2-nk}\\
&=\frac{q^{n(n+1)}{(-1)^n}}{[2n+1]_q!}\frac12
\prod_{j=0}^{2n}(1+q^jz)\Big\vert_{z=q^{-n}}\\
&=\frac{q^{\binom{n+1}2}{(-1)^n}}{[2n+1]_q!}
\prod_{i=1}^n(1+q^i)^2.
\end{align*}

Although we do not need it, we also mention the dual formula
\begin{align*}
[z^{2n+1}]\sin_{q}z\cos_{1/q}z
=\frac{q^{\binom{n}2}{(-1)^n}}{[2n+1]_q!}
\prod_{i=1}^n(1+q^i)^2.
\end{align*}

A similar computation gives the result ($n\ge1$)

\begin{align*}
[z^{2n}]\cos_{q}z\cos_{1/q}z=-[z^{2n}]\sin_{q}z\sin_{1/q}z
=\frac{q^{\binom{n}2}{(-1)^n}}{[2n]_q!}
\prod_{i=1}^{n-1}(1+q^i)^2(1+q^n).
\end{align*}

Now we write $\tan_qz=\frac{\sin_{q}z\cos_{1/q}z}{\cos_{q}z\cos_{1/q}z}$ and 
thus
\begin{equation*}
{\cos_{q}z\cos_{1/q}z}\sum_{n\ge0}\frac{T_{2n+1}(q)}{[2n+1]_q!}
z^{2n+1}={\sin_{q}z\cos_{1/q}z}.
\end{equation*}

Comparing coefficients, we find

\begin{align}\begin{split}\label{tanrec}
T_{2n+1}(q)+\sum_{k=1}^{n}\gauss{2n+1}{2k}q^{\binom k2}(-1)^k
\prod_{i=1}^{k-1}(1+q^i)^2(1+q^k)T_{2n+1-2k}(q)
\\=q^{\binom n2}(-1)^n\prod_{i=1}^{n}(1+q^i)^2.
\end{split}\end{align}
The induction argument is as in \cite{AnGe78};
$T_{2n+1-2k}(q)$ has a factor $\prod_{i=1}^{n-k}(1+q^i)$ and, according again
to \cite{AnGe78},
\begin{equation*}
\gauss{2n+1}{2k}\frac{\prod_{i=1}^{k}(1+q^i)}{\prod_{i=n-k+1}^{n}(1+q^i)}
\end{equation*}
is still a polynomial.
The two factors $\prod_{i=n-k+1}^{n}(1+q^i)$ and $\prod_{i=1}^{n-k}(1+q^i)$ 
mean that everything in (\ref{tanrec}) must be divisible by
$\prod_{i=1}^{n}(1+q^i)$, and this finishes the proof.

It is likely that
stronger results as in \cite{Foata81} hold, but we have not investigated 
that. 
\end{proof}

The new $q$--secant numbers do not enjoy any divisibility results that are
worthwhile to report; for the classical ones, see \cite{AnFo80}.

{\bf Remark.} The paper \cite{IsRaSt99} has a $q$--exponential function
\begin{equation*}
\mathcal{E}_q:=\sum_{n\ge0}\frac{q^{n^2/4}z^n}{(q;q)_n}.
\end{equation*}

Plugging in $iz(1-q)$ for $z$ and taking real parts would result in the $q$--cosine
with factor $q^{n^2}$. To get the corresponding $q$--sine,
replace $z$ by $izq(1-q)$, take the imaginary part and multiply by $q^{1/4}$.
We consider that merely to be a curiosity, not being of much help.

{\bf Acknowledgment.} I want to thank Dominique Foata for several
pointers to the literature. 

\bibliographystyle{plain}


\begin{thebibliography}{10}

\bibitem{Andrews76}
G.~Andrews.
\newblock {\em The Theory of Partitions}, volume~2 of {\em Encyclopedia of
  Mathematics and its Applications}.
\newblock Addison--Wesley, 1976.

\bibitem{AnAsRo99}
G.~Andrews, R.~Askey, and R.~Roy.
\newblock {\em Special Functions}, volume~71 of {\em Encyclopedia of
  Mathematics and its Applications}.
\newblock Cambridge University Press, 1999.

\bibitem{AnFo80}
G.~Andrews and D.~Foata.
\newblock Congruences for the $q$--secant numbers.
\newblock {\em European Journal of Combinatorics}, 1:283--287, 1980.

\bibitem{AnGe78}
G.~Andrews and I.~Gessel.
\newblock Divisibility properties of the $q$--tangent numbers.
\newblock {\em Proceedings of the American Mathematical Society}, 68:380--384,
  1978.

\bibitem{Flajolet80}
P.~Flajolet.
\newblock Combinatorial aspects of continued fractions.
\newblock {\em Discrete Applied Mathematics}, 39:207--229, 1992.

\bibitem{FlPr87}
P.~Flajolet and H.~{P}rodinger.
\newblock Level number sequences for trees.
\newblock {\em Discrete Mathematics}, 65:149--156, 1987.

\bibitem{Foata81}
D.~Foata.
\newblock Further divisibility properties of the $q$--tangent numbers.
\newblock {\em Proceedings of the American Mathematical Society}, 81:143--148,
  1981.

\bibitem{GrKnPa94}
R.~L. Graham, D.~E. Knuth, and O.~Patashnik.
\newblock {\em Concrete Mathematics (Second Edition)}.
\newblock Addison Wesley, 1994.

\bibitem{HaRaZe99}
G.~N. Han, A.~Randrianarivony, and J.~Zeng.
\newblock Un autre $q$--analogue des nombres d'{E}uler.
\newblock {\em 42e S\'eminaire Lotharingien}, [B42e]:22 pages, 1999.

\bibitem{IsRaSt99}
M.~Ismail, M.~Rahman, and D.~Stanton.
\newblock Quadratic $q$--exponentials and connection coefficient problems.
\newblock {\em preprint}, 12 pages, 1999.

\bibitem{Jackson04}
F.~H. Jackson.
\newblock A basic--sine and cosine with symbolic solutions of certain
  differential equations.
\newblock {\em Proc. Edinburg Math. Soc.}, 22:28--39, 1904.

\bibitem{KnPr98}
A.~Knopfmacher and H.~{P}rodinger.
\newblock On {C}arlitz compositions.
\newblock {\em European Journal on Combinatorics}, 19:579--589, 1998.

\end{thebibliography}
\end{document}